\documentclass[12pt,leqno]{amsart}
\usepackage[dvips]{graphicx}
\usepackage{amsfonts}
\usepackage{amssymb}
\usepackage{amsmath}
\usepackage{amscd}
\usepackage{graphicx}

\voffset-1cm \markleft{\rm CNAM Probabilit\'es. Chapitre 2. G\'en\'eralit\'es } 
\newtheorem{theorem}{Theorem}

\newtheorem{lemma}[theorem]{Lemma}

\newtheorem{proposition}[theorem]{Proposition}

\newcommand\R{\mathbb{R}}

\newcommand\K{\mathbb{K}}
\newcommand\F{\mathbb{F}}

\newcommand\ds{\displaystyle}
\newcommand\pf{\noindent{\it Proof. }}
\title{2-dimensional algebras. Applications to Jordan, $G$-associative and Hom-associative algebras}

\author{Elisabeth Remm}
\email{Elisabeth.Remm@uha.fr}
\address{Laboratoire de Math\'ematiques et Applications,
        Universit\'e de Haute Alsace,\\ Facult\'e des Sciences et
        Techniques, 4, rue des Fr\`eres Lumi\`ere,\\
        68093~Mulhouse~cedex, France.}
\author{Michel Goze} 
\email{goze.rac@gmail.com}
\address{Ramm Algebra Center, 4 rue de Cluny, F.68800 Rammersmatt}


\begin{document}

\begin{abstract}
We classify, up to isomorphism, the $2$-dimensional algebras over a field $\K$. We focuse  also on the case of characteristic $2$,  identifying the matrices of $GL(2,\F_2)$ with the elements of the symmetric group $\Sigma_3$. The classification  is then given by the study of the orbits of this group on a $3$-dimensional plane, viewed as a Fano plane. As applications, we establish  classifications of Jordan algebras, algebras of Lie type or Hom-Associative algebras.
  \end{abstract}




\bibliographystyle{plain}

\maketitle

\section{Introduction}

An algebra $\mathcal{A}$ over a field $\K$ is $\K$-vector space  equipped with a product  which corresponds to a bilinear map on $\mathcal{A}$ with values in $\mathcal{A}.$  For a given dimension, one of the basic problems is the determination  up to linear isomorphism of all these algebras. Sub classes of algebras where widely studied . These subclasses where often obtained setting a quadratic relation on $\mu.$  Among other examples of such classes are Lie algebras (in this case  $\mu$ is skewsymmetric and satisfies Jacodi identity), associative algebras, Lie-admissibles algebras, Pre-Lie algebras in particular. In all these examples, classifications where established in a general frame work, that is, with no other hypothesis on these classes  and only in very small dimensions. For example for Lie algebras, we know the general classifications up to the dimension $6$. In bigger dimension we impose additional algebraic properties if we hope to continue this classification.  For example simple Lie algebras are fully classified since the work of Killing and Cartan, in any dimension. Unfortunately it is more and less the only solved case. If we consider complexe nilpotent Lie algebras, the classification is known only up to the dimension $7$.  It is the same for the associatives algebras. If we are only interested in general algebras, the only known cases are the dimension $2$ and $3.$ It is true that the problem is equivalent to  the classification of tensors of type $(2,1)$ on a finite dimensional vector space. We are then facing to a basic multilinear algebra problem which is subject to a lack of informations on the tensors.

Here we reconsider this problem from the beginning, that is in dimension $2$. This work is certainly not the first one of the subject. There is for example the work of Petersson. Our approach is not similar. We are not fully  interested by the classification up to isomorphism but by the determination of subclasses, minimal in a certain sense, which are invariant up to isomorphism. The motivation  comes from the constatation of what happen in greater dimensions for nilpotent Lie algebras for example In this case, the classification is established in dimension $7$ but  quasi unusable in its present forme. This means that if we have  a precise example of nilpotent Lie algebra of this dimension, it is long and fastidious to recognize it in the given list because most of the time it is not adapted to the invariants used to established the classification. Moreover the length of the list can be  puzzling. In greater dimensions, the number of isomorphy classes, the need to write invariant parametrized families seems to be an unrealistic goal.  Hence the idea to reduce the classification problem to a determination of invariant classes. This is the aim of this work.  However we will established the link with Petersson's work. Our approach is quite basic. In characteristic different from $2$, we decompose a tensor   $\mu$ as a skewsymmetric and symmetric one. Since the skewsymmetric case  is elementary, we classify those which are symmetric  modulo the automorphism group of the associated skeysymmetric law. In characteristic  $2$, the problem is equivalent to the determination of the orbits of the Fano plane modulo the symmetric group. Finally, we use these results to describe or find again certain classes of algebras whose a direct approach is rather difficult. In particular, we determine the $2$-dimensional  Jordan algebras and we find again the results of \cite{Ancochea}, the $G$-associative algebras  and the Hom-associative algebras.

We have begun the study of the determination of general algebras in \cite{GRMaroc} which was  specially an introduction to a more precise work developed in this paper but with the same idea to describe "minimal" families invariant by isomorphism rather than a precise list for which the use is difficult. Recently, we were acquainted with the work of Pertersson, based on an Kaplansky result which permits to describe all the algebras from some unital algebras and to give isomorphism criteria. We try in this paper to look our description in a Petersson point of view. We note also a recent work, on the same subject of H. Ahmed, U. Bekbaev and I. Rakhimov \cite{Be}.  

\medskip

\section{Generalities}
Let $\K$ be a field whose characteristic will be precise later. 
An algebra over a field $\K$ is a $\K$-vector space $V$ with a multiplication given by a bilinear map
$$\mu: V \times V\rightarrow  V.$$
We denote by $A=(V,\mu)$ a $\K$-algebra structure on $V$ with multiplication $\mu$.
Throughout this paper we fix the vector space $V$. Since we are interested by the $2$-dimensional case we could assume that $V=\K^2$. 
Two $\K$-algebras $A=(V,\mu )$ and $A'=(V,\mu ')$
are  isomorphic if there
is a linear isomorphism
$$f : V\rightarrow  V$$
such as
$$f(\mu (X,Y))=\mu '(f(X),f(Y))$$
for all $X,Y \in V.$ The classification of $2$-dimensional $\K$-algebras is then equivalent to the classification of
bilinear maps on $V=\K^2$ with values in $V$. Let $\{e_1,e_2\}$ be a fixed basis of $V$. A  general bilinear map $\mu$ has the following expression
$$
\left\{
\begin{array}{l}
\mu(e_1,e_1)=\alpha_1e_1+\beta _1e_2, \\
\mu(e_1,e_2)=\alpha_2e_1+\beta _2e_2, \\
\mu(e_2,e_1)=\alpha_3e_1+\beta _3e_2, \\
 \mu(e_2,e_2)=\alpha_4e_1+\beta _4e_2, \\
 \end{array}
 \right.
 $$
 and it is defined by $8$ parameters. Let $f$ be a linear isomorphism of $V$. In the given basis, its matrix $M$ is non degenerate. If we put
 $$M=\left(
 \begin{array}{rr}
 a & b\\
 c & d
 \end{array}
 \right)
 $$
 then
 $$M^{-1}=
 \frac{1}{\Delta }
 \left(\begin{array}{rr}
 d & -b\\
 -c & a
 \end{array}
 \right)
 $$
 with $\Delta =ad-bc \not= 0.$ The isomorphic multiplication
 $$\mu'=f^{-1}\circ \mu \circ (f \times  f)$$
 satisfies
 $$
\left\{
\begin{array}{l}
\mu'(e_1,e_1)=\alpha'_1e_1+\beta' _1e_2, \\
\mu'(e_1,e_2)=\alpha'_2e_1+\beta' _2e_2, \\
\mu'(e_2,e_1)=\alpha'_3e_1+\beta' _3e_2, \\
 \mu'(e_2,e_2)=\alpha'_4e_1+\beta' _4e_2, \\
 \end{array}
 \right.
 $$
 with
 \begin{equation}\label{change}
\left\{
\begin{array}{l}
\alpha'_1= (a^2\alpha _1+ac\alpha _2+ac\alpha _3+c^2\alpha _4)\frac{d}{\Delta }
-(a^2\beta  _1+ac\beta  _2+ac\beta  _3+c^2\beta  _4)\frac{b}{\Delta }\\
\beta' _1= -(a^2\alpha _1+ac\alpha _2+ac\alpha _3+c^2\alpha _4)\frac{c}{\Delta }
+(a^2\beta  _1+ac\beta  _2+ac\beta  _3+c^2\beta  _4)\frac{a}{\Delta }\\
\alpha'_2=(ab\alpha _1+ad\alpha _2+bc\alpha _3+cd\alpha _4)\frac{d}{\Delta }
-(ab\beta  _1+ad\beta  _2+bc\beta  _3+cd\beta  _4)\frac{b}{\Delta }\\
\beta' _2= -(ab\alpha _1+ad\alpha _2+bc\alpha _3+cd\alpha _4)\frac{c}{\Delta }
+(ab\beta  _1+ad\beta  _2+bc\beta  _3+cd\beta  _4)\frac{a}{\Delta }\\
\alpha'_3=(ab\alpha _1+bc\alpha _2+ad\alpha _3+cd\alpha _4)\frac{d}{\Delta }
-(ab\beta  _1+bc\beta  _2+ad\beta  _3+cd\beta  _4)\frac{b}{\Delta }\\
\beta' _3= -(ab\alpha _1+bc\alpha _2+ad\alpha _3+cd\alpha _4)\frac{c}{\Delta }
+(ab\beta  _1+bc\beta  _2+ad\beta  _3+cd\beta  _4)\frac{a}{\Delta }\\
\alpha'_4=(b^2\alpha _1+bd\alpha _2+bd\alpha _3+d^2\alpha _4)\frac{d}{\Delta }
-(b^2\beta  _1+bd\beta  _2+bd\beta  _3+d^2\beta  _4)\frac{b}{\Delta }\\
\beta' _4= -(b^2\alpha _1+bd\alpha _2+bd\alpha _3+d^2\alpha _4)\frac{c}{\Delta }
+(b^2\beta  _1+bd\beta  _2+bd\beta  _3+d^2\beta  _4)\frac{a}{\Delta }\\
 \end{array}
 \right.
\end{equation}
 These formulae describe an action of the linear group $GL(2,\K)$ on $\K^8$ parameterized by the structure constants $(\alpha _i,\beta _i)$,
 $i=1,2,3,4$ and the problem of classification consists in describing an element of each orbit.

\section{Algebras over a field of characteristic different from $2$}

We assume in this section that $char(\K) \neq 2$. We consider the bilinear map $\mu_a$ and $\mu_s$ given by
$$\mu_a(X,Y)=\displaystyle \frac{\mu(X,Y)-\mu(Y,X)}{2}, \ \  \mu_s(X,Y)=\displaystyle \frac{\mu(X,Y)+\mu(Y,X)}{2}
$$
for all $X,Y \in V$. The multiplication  $\mu_a$ is skew-symmetric and it is a Lie multiplication (any skew-symmetric bilinear application in $\K^2$ is a Lie bracket). It is isomorphic to one of the following

\begin{enumerate}
  \item $\mu_a^1(e_1,e_2)=e_1,$
  \item $\mu_a^2=0.$
\end{enumerate}

In fact, if $\mu_a$ is not trivial, thus $\mu_a(e_1,e_2)=\alpha e_1+\beta e_2$. If $\alpha \neq 0$, we consider the change of basis
$$e'_1=\alpha e_1+\beta e_2, \ e'_2=\alpha^{-1}e_2.$$
We have $\mu_a(e'_1,e'_2)=\mu_a(\alpha e_1+\beta e_2,\alpha^{-1}e_2)=\mu_a(e_1,e_2)=\alpha e_1+\beta e_2=e'_1.$

\noindent If $\alpha =0$, thus $\beta  \neq 0$ and we take $$e'_1=e_2, \ e'_2=-\beta ^{-1}e_1.$$ This gives $\mu_a(e'_1,e'_2)=\mu_a(e_2,-\beta ^{-1}e_1)=\beta ^{-1}\beta e_2=e_2=e'_1.$ In any case, if $\mu_a \neq 0$, then it is isomorphic to $\mu_a^1.$

\subsection{Case $\mu^1_a(e_1,e_2)=e_1$}

\medskip

An automorphism of the Lie algebra $(A,\mu_a^1)$ is a linear isomorphism $f \in GL(2,\K)$ such that
$$f(\mu_a^1(X,Y))=\mu_a^1(f(X),f(Y))$$
for every $X,Y \in A$. The set of automorphisms of this Lie algebra is  denoted by $Aut(\mu_a^1)$.
\begin{lemma}
We have
$$ Aut(\mu_a^1)=\{M=\left(
 \begin{array}{ll}
 a & b\\
 0& 1
 \end{array}
 \right), \ a,b\in \K, \ a\neq 0\}.
 $$
 \end{lemma}
\pf  In fact, assume that $M=\left(
 \begin{array}{ll}
 a & b\\
 c& d
 \end{array}
 \right)$ is the matrix of the automorphism $f$ in the given basis $\{e_1,e_2\}$.  Then
 $$f(\mu_a^1(e_1,e_2))=f(e_1)=ae_1+ce_2,$$
 and
 $$\mu_a^1(f(e_1),f(e_2))=\mu_a^1(ae_1+ce_2,be_1+de_2)=(ad-bc)e_1.$$
 Then
 $$c=0, \ a=ad.$$
 But $\det M= ad \neq 0$ so $a=ad$ implies that $d=1.$ This gives the lemma. $\Box$

\medskip

Let $\mu$ be a general multiplication of $2$-dimensional $\K$-algebra such that $\mu_a$ is isomorphic to $\mu_a^1$. It is isomorphic to a the bilinear map (always denoted by $\mu$) whose structural constants are given by
$$
\left\{
\begin{array}{l}
\mu(e_1,e_1)=\alpha_1e_1+\beta _1e_2, \\
\mu(e_1,e_2)=(\alpha_2+1)e_1+\beta _2e_2, \\
\mu(e_2,e_1)=(\alpha_2-1)e_1+\beta _2e_2, \\
 \mu(e_2,e_2)=\alpha_4e_1+\beta _4e_2.\\
 \end{array}
 \right.
 $$
The classification, up to isomorphism, of the Lie algebras $(V,\mu)$ such that $\mu_a$ is isomorphic to $\mu_a^1$ is equivalent to the classification up an isomorphism belonging to $Aut(\mu_a^1)$ of the abelian algebras isomorphic to
 $$
\left\{
\begin{array}{l}
\mu_s(e_1,e_1)=\alpha_1e_1+\beta _1e_2, \\
\mu_s(e_1,e_2)=\mu_s(e_2,e_1)=\alpha_2e_1+\beta _2e_2, \\
 \mu_s(e_2,e_2)=\alpha_4e_1+\beta _4e_2, \\
 \end{array}
 \right.
 $$
 In this case (\ref{change}) is reduced to
  \begin{equation}\label{c1}
\left\{
\begin{array}{l}
\alpha'_1= a\alpha _1-ab\beta  _1,\\
\beta' _1=a^2\beta  _1,\\
\alpha'_2=\alpha'_3=b\alpha _1+\alpha _2-b^2\beta  _1-b\beta  _2,\\
\beta' _2= \beta'_3=ab\beta  _1+a\beta  _2,\\
\alpha'_4=\ds (b^2\alpha _1+2b\alpha _2+\alpha _4
-b^3\beta  _1-2b^2\beta  _2-b\beta  _4)\frac{1}{a},\\
\beta' _4=
b^2\beta  _1+2b\beta  _2+\beta  _4.\\
 \end{array}
 \right.
\end{equation}
\begin{enumerate}
  \item Assume that $\beta_1 \neq 0.$

  $\bullet$ Suppose that $\K$ is algebraically closed and consider the isomorphism $\begin{pmatrix}
\sqrt{\beta_1}      & \frac{\alpha_1}{\beta_1}   \\
 0     & 1
\end{pmatrix}.$
The isomorphic algebra is such that
$\alpha'_1= 0$ and $\beta' _1=1.$
We deduce that in this case $\mu_s$ is isomorphic to
 $$
\left\{
\begin{array}{l}
\mu_s(e_1,e_1)=e_2, \\
\mu_s(e_1,e_2)=\mu_s(e_2,e_1)=\alpha_2e_1+\beta _2e_2, \\
 \mu_s(e_2,e_2)=\alpha_4e_1+\beta _4e_2. \\
 \end{array}
 \right.
 $$
  Then $\mu $ is isomorphic to
  $$
\left\{
\begin{array}{l}
\mu^1_{\alpha_2,\beta_2,\alpha_4,\beta_4}(e_1,e_1)=e_2, \\
\mu^1_{\alpha_2,\beta_2,\alpha_4,\beta_4}(e_1,e_2)=(\alpha_2+1)e_1+\beta _2e_2, \\
\mu^1_{\alpha_2,\beta_2,\alpha_4,\beta_4}(e_2,e_1)=(\alpha_2-1)e_1+\beta _2e_2, \\
 \mu^1_{\alpha_2,\beta_2,\alpha_4,\beta_4}(e_2,e_2)=\alpha_4e_1+\beta _4e_2, \\
 \end{array}
 \right.
 $$
 with $\alpha_2,\beta_2,\alpha_4,\beta_4 \in \K.$

 \medskip

 $\bullet$ If $\K$ is not algebraically closed (for example if $\K$ is a finite field), let $\K^{*2} $ be the multiplicative subgroup of elements $a^2$ with $a \in \K$. In this case $\mu$ is isomorphic to a Lie bracket belonging to the $4$ parameters family:
$$
\left\{
\begin{array}{l}
\varphi^{1,\lambda}_{\alpha_2,\beta_2,\alpha_4,\beta_4}(e_1,e_1)=\lambda e_2, \\
\varphi^{1,\lambda}_{\alpha_2,\beta_2,\alpha_4,\beta_4}(e_1,e_2)=(\alpha_2+1)e_1+\beta _2e_2, \\
\varphi^{1,\lambda}_{\alpha_2,\beta_2,\alpha_4,\beta_4}(e_2,e_1)=(\alpha_2-1)e_1+\beta _2e_2, \\
\varphi^{1,\lambda}_{\alpha_2,\beta_2,\alpha_4,\beta_4}(e_2,e_2)=\alpha_4e_1+\beta _4e_2, \\
 \end{array}
 \right.
 $$
 with $\alpha_2,\beta_2,\alpha_4,\beta_4 \in \K$ and $\lambda \in \K/ \K^{*^2} $. For example, if $\K=\R$, then $\lambda \in \{-1,1\}$.

  \medskip

  \item Assume $\beta_1=0,$ $\beta_2 \neq 0.$ In this case (\ref{change}) is reduced to
  \begin{equation}\label{change2}
\left\{
\begin{array}{l}
\alpha'_1= a\alpha _1,\\
\beta'_1=0,\\
\alpha'_2=b\alpha _1+\alpha _2
-b\beta  _2,\\
\beta' _2=
a\beta  _2,\\
\alpha'_4=\ds (b^2\alpha _1+2b\alpha _2+\alpha _4
-2b^2\beta  _2-b\beta  _4)\frac{1}{a},\\
\beta' _4=
2b\beta  _2+\beta  _4.\\
 \end{array}
 \right.
\end{equation}
and taking $b=-\beta_4/2\beta_2$ and $a=\beta_2^{-1}$, we see that $\mu_s$ is isomorphic to
  $$
\left\{
\begin{array}{l}
\mu_s(e_1,e_1)=\alpha _1e_1, \\
\mu_s(e_1,e_2)=\mu_s(e_2,e_1)=\alpha_2e_1+e_2, \\
 \mu_s(e_2,e_2)=\alpha_4e_1.\\
 \end{array}
 \right.
 $$
 We obtain the following multiplication,  $\K$ being algebraically closed or not:
 $$
\left\{
\begin{array}{l}
\mu^2_{\alpha_1,\alpha_2,\alpha_4}(e_1,e_1)=\alpha _1e_1, \\
\mu^2_{\alpha_1,\alpha_2,\alpha_4}(e_1,e_2)=(\alpha_2+1)e_1+e_2, \\
\mu^2_{\alpha_1,\alpha_2,\alpha_4}(e_2,e_1)=(\alpha_2-1)e_1+e_2, \\
 \mu^2_{\alpha_1,\alpha_2,\alpha_4}(e_2,e_2)=\alpha_4e_1. \\
 \end{array}
 \right.
 $$
  \item Assume now that $\beta_1=\beta_2=0,\alpha_1 \neq 0$. In this case (\ref{change}) is reduced to
  \begin{equation}\label{change3}
\left\{
\begin{array}{l}
\alpha'_1= a\alpha _1\\
\beta' _1=\beta' _2=0\\
\alpha'_2=b\alpha _1+\alpha _2\\
\alpha'_4=\ds (b^2\alpha _1+2b\alpha _2+\alpha _4
-b\beta  _4)\frac{1}{a}\\
\beta' _4=
\beta  _4.\\
 \end{array}
 \right.
\end{equation}
and taking $b=-\alpha_2/\alpha_1$ and $a=\alpha_1^{-1}$,  we obtain $\alpha'_2=0$ and $\alpha'_1=1.$ In this case, $\mu$ is isomorphic to
$$
\left\{
\begin{array}{l}
\mu^3_{\alpha_4,\beta_4}(e_1,e_1)=e_1, \\
\mu^3_{\alpha_4,\beta_4}(e_1,e_2)=e_1, \\
\mu^3_{\alpha_4,\beta_4}(e_2,e_1)=-e_1, \\
 \mu^3_{\alpha_4,\beta_4}(e_2,e_2)=\alpha_4e_1+\beta _4e_2. \\
 \end{array}
 \right.
 $$

  \item
Assume now that $\beta_1=\beta_2=0,\alpha_1=0,2\alpha_2 -\beta_4\neq 0$. In this case, considering $b=-\alpha_4/(2\alpha_2-\beta_4)$, the Lie bracket $\mu$ is isomorphic to
$$
\left\{
\begin{array}{l}
\mu^4_{\alpha_2,\beta_4}(e_1,e_1)=0, \\
\mu^4_{\alpha_2,\beta_4}(e_1,e_2)=(\alpha_2+1)e_1, \\
\mu^4_{\alpha_2,\beta_4}(e_2,e_1)=(\alpha_2-1)e_1,\\
 \mu^4_{\alpha_2,\beta_4}(e_2,e_2)=\beta _4e_2, \\
 \end{array}
 \right.
 $$

  \item Assume now that $\beta_1=\beta_2=0,\alpha_1=0,2\alpha_2 -\beta_4= 0, \alpha_4\neq 0$. The Lie bracket $\mu$ is isomorphic to
$$
\left\{
\begin{array}{l}
\mu^5_{\alpha_2}(e_1,e_1)=0, \\
\mu^5_{\alpha_2}(e_1,e_2)=(\alpha_2+1)e_1, \\
\mu^5_{\alpha_2}(e_2,e_1)=(\alpha_2-1)e_1, \\
 \mu^5_{\alpha_2}(e_2,e_2)= e_1+2\alpha_2e_2, \\
 \end{array}
 \right.
 $$

  \item If $\beta_1=\beta_2=0,\alpha_1=0,2\alpha_2 -\beta_4= 0, \alpha_4=0$, then $\mu$ is isomorphic to $\mu^4_{\alpha_2,\beta_4}$ with
  $\beta_4=2\alpha_2$
\end{enumerate}

\begin{theorem} \label{nc}
Any 2-dimensional non commutative algebras isomorphic to one of the following algebras:
\begin{itemize}
  \item If $\K$ is algebraically closed $$
\left\{
\begin{array}{l}
\mu^1_{\alpha_2,\beta_2,\alpha_4,\beta_4}(e_1,e_1)=e_2, \\
\mu^1_{\alpha_2,\beta_2,\alpha_4,\beta_4}(e_1,e_2)=(\alpha_2+1)e_1+\beta _2e_2, \\
\mu^1_{\alpha_2,\beta_2,\alpha_4,\beta_4}(e_2,e_1)=(\alpha_2-1)e_1+\beta _2e_2, \\
 \mu^1_{\alpha_2,\beta_2,\alpha_4,\beta_4}(e_2,e_2)=\alpha_4e_1+\beta _4e_2. \\
 \end{array}
 \right.
\qquad
\left\{
\begin{array}{l}
\mu^2_{\alpha_1,\alpha_2,\alpha_4}(e_1,e_1)=\alpha _1e_1, \\
\mu^2_{\alpha_1,\alpha_2,\alpha_4}(e_1,e_2)=(\alpha_2+1)e_1+e_2, \\
\mu^2_{\alpha_1,\alpha_2,\alpha_4}(e_2,e_1)=(\alpha_2-1)e_1+e_2, \\
 \mu^2_{\alpha_1,\alpha_2,\alpha_4}(e_2,e_2)=\alpha_4e_1. \\
 \end{array}
 \right.
 $$
 $$
\left\{
\begin{array}{l}
\mu^3_{\alpha_4,\beta_4}(e_1,e_1)=e_1, \\
\mu^3_{\alpha_4,\beta_4}(e_1,e_2)=e_1, \\
\mu^3_{\alpha_4,\beta_4}(e_2,e_1)=-e_1, \\
 \mu^3_{\alpha_4,\beta_4}(e_2,e_2)=\alpha_4e_1+\beta _4e_2. \\
 \end{array}
 \right.
 \left\{
\begin{array}{l}
\mu^4_{\alpha_2,\beta_4}(e_1,e_1)=0, \\
\mu^4_{\alpha_2,\beta_4}(e_1,e_2)=(\alpha_2+1)e_1, \\
\mu^4_{\alpha_2,\beta_4}(e_2,e_1)=(\alpha_2-1)e_1, \\
 \mu^4_{\alpha_2,\beta_4}(e_2,e_2)=\beta _4e_2, \\
 \end{array}
 \right.
\left\{
\begin{array}{l}
\mu^5_{\alpha_2}(e_1,e_1)=0, \\
\mu^5_{\alpha_2}(e_1,e_2)=(\alpha_2+1)e_1, \\
\mu^5_{\alpha_2}(e_2,e_1)=(\alpha_2-1)e_1, \\
 \mu^5_{\alpha_2}(e_2,e_2)= e_1+2\alpha_2e_2. \\
 \end{array}
 \right.
 $$
with $\alpha_i,\beta_i \in \K$.
  \item If $\K$ is not algebraically closed
   $$
\left\{
\begin{array}{l}
\varphi^{1,\lambda}_{\alpha_2,\beta_2,\alpha_4,\beta_4}(e_1,e_1)=\lambda e_2, \\
\varphi^{1,\lambda}_{\alpha_2,\beta_2,\alpha_4,\beta_4}(e_1,e_2)=(\alpha_2+1)e_1+\beta _2e_2, \\
\varphi^{1,\lambda}_{\alpha_2,\beta_2,\alpha_4,\beta_4}(e_2,e_1)=(\alpha_2-1)e_1+\beta _2e_2, \\
\varphi^{1,\lambda}_{\alpha_2,\beta_2,\alpha_4,\beta_4}(e_2,e_2)=\alpha_4e_1+\beta _4e_2, \\
 \end{array}
 \right., \   \mu^2_{\alpha_1,\alpha_2,\alpha_4} \ ,  \quad \mu^3_{\alpha_4,\beta_4} \ ,  \quad \mu^4_{\alpha_2,\beta_4} \ ,  \quad \mu^5_{\alpha_2} \ ,$$
 $\alpha_i,\beta_i \in \K, \lambda \in \K/(\K^*)^2.$
\end{itemize}
\end{theorem}

\medskip

Let us make the link with the results of Petersson (\cite{Pet}). The main idea of this work is to construct algebras from unital algebra. 
Recall that an algebra $A=(V,\mu)$ is called unital if there exists $1 \in V$ such that $\mu(1,X)=\mu(X,1)=X$ for any $X \in V$.
\begin{lemma}
If $\mu_a$ is not trivial, then $A$ is not unital.
\end{lemma}
\pf Assume that there exists $1$ satisfying $\mu(1,X)=\mu(X,1)=X$, then
$$0=\mu(1,X)-\mu(X,1)=\mu_a(1,X)-\mu_a(X,1)=2\mu_a(1,X)$$
for any $X \in V$. Then $\mu_a(1,X)=0$ for any $X$ and $1$ is in the center of $A_a=(V,\mu_a)$. But if $\mu_a$ is not trivial, the center of $A_a$ is reduce to $\{0\}$. The algebra $A$ cannot be unital.

The algebra $A=(V,\mu)$ is called regular if there exists $U,T\in V$ such that the linear applications
$$L_U:X \rightarrow \mu(U,X), \ \ R_T: X \rightarrow \mu(X,T)$$
are linear isomorphisms. From \cite{Pet}, for any regular algebra $A=(V,\mu)$ there exist a unique, up an isomorphism, unital algebra $B=(V,\mu_u)$ and two linear isomorphisms $f, g$ of $V$ such that
$$\mu(X,Y)= \mu_u(f(X),g(Y)$$
for any $X,Y \in V$. The algebra $B$ is called the unital heart of $A$. To compare Theorem \ref{nc} with the Petersson results, we have to determine the regular algebras. Let us consider the first family. The application $L_U$ is not regular for any $U$ if and only if its determinant is identically null that is
$$\alpha_2=-1,  \ \alpha_4=-2\beta_2, \ \beta_4=\beta_2^2.$$
Likewise $R_T$ is not regular for any $T$ if and only if its determinant is identically null that is
$$\alpha_2=1,  \ \alpha_4=2\beta_2, \ \beta_4=\beta_2^2.$$

We deduce that any algebra $A^1_{\alpha_2,\beta_2,\alpha_4,\beta_4}=(V,
\mu^1_{\alpha_2,\beta_2,\alpha_4,\beta_4})$ is regular except the algebras given by
$$
\left\{
\begin{array}{l}
\mu^1_{-1,\beta_2,-2\beta_2,\beta_2^2}(e_1,e_1)=e_2, \\
\mu^1_{-1,\beta_2,-2\beta_2,\beta_2^2}(e_1,e_2)=\beta _2e_2, \\
\mu^1_{-1,\beta_2,-2\beta_2,\beta_2^2}(e_2,e_1)=-2e_1+\beta _2e_2, \\
 \mu^1_{-1,\beta_2,-2\beta_2,\beta_2^2}(e_2,e_2)=-2\beta_2e_1+\beta _2^2e_2. \\
 \end{array}
 \right.
 \ \
 \left\{
\begin{array}{l}
\mu^1_{1,\beta_2,2\beta_2,\beta_2^2}(e_1,e_1)=e_2, \\
\mu^1_{1,\beta_2,2\beta_2,\beta_2^2}(e_1,e_2)=2e_1+\beta _2e_2, \\
\mu^1_{1,\beta_2,2\beta_2,\beta_2^2}(e_2,e_1)=\beta _2e_2, \\
 \mu^1_{1,\beta_2,2\beta_2,\beta_2^2}(e_2,e_2)=2\beta_2e_1+\beta _2^2e_2. \\
 \end{array}
 \right.
$$
Let us note that $A^1_{-1,\beta_2,-2\beta_2,\beta_2^2}$ is left-singular but right-regular and $A^1_{1,\beta_2,2\beta_2,\beta_2^2}$
is right-singular and left-regular. An algebra which is left and right singular is called bi-singular. We can summarize the results in the following array:
\begin{enumerate}
  \item $A^1_{\alpha_2,\beta_2,\alpha_4,\beta_4}$ regular except $A^1_{-1,\beta_2,-2\beta_2,\beta_2^2}$ and $A^1_{1,\beta_2,2\beta_2,\beta_2^2}$.\\
  \item $A^1_{-1,\beta_2,-2\beta_2,\beta_2^2}$ is left-singular and right-regular \\
  \item $A^1_{1,\beta_2,2\beta_2,\beta_2^2}$ is right-singular and left-regular,\\
  \item $A^2_{\alpha_1,\alpha_2,\alpha_4}$ is regular,\\
  \item $A^3_{\alpha_4,\beta_4}$ is regular except $A^3_{\alpha_4,0}$, \\
  \item $A^3_{\alpha_4,0}$ is bisingular. \\
   \item $A^4_{\alpha_2,\beta_4}$ is regular except $A^4_{\alpha_2,0}$, $A^4_{1,\beta_4}$, $A^4_{-1,\beta_4}$ \\
\item $A^4_{\alpha_2,0}$ is bisingular,
\item $A^4_{1,\beta_4}$ is left-singular and right-regular as soon as $\beta_4 \neq 0$,
\item $A^4_{-1,\beta_4}$ is left-regular and right-singular as soon as $\beta_4 \neq 0$,
 \item $A^5_{\alpha_2}$ is regular except for $\alpha_2=0, 1$ or $-1$, \\
\item $A^5_{0}$ is bisingular,
\item $A^5_{1}$ is left-singular and right-regular as soon as $\beta_4 \neq 0$,
\item $A^5_{-1}$ is left-regular and right-singular as soon as $\beta_4 \neq 0$,
\end{enumerate}

We deduce
\begin{proposition}
We consider the following algebras
\begin{enumerate}
  \item $A^1_{\alpha_2,\beta_2,\alpha_4,\beta_4}$ with $(\alpha_2,\beta_2,\alpha_4,\beta_4) \neq (-1,\beta_2,-2\beta_2,\beta_2^2)$ or $(1,\beta_2,2\beta_2,\beta_2^2),$
  \item $A^2_{\alpha_1,\alpha_2,\alpha_4}$,
  \item $A^4_{\alpha_2,\beta_4}$ with $(\alpha_2,\beta_4) \neq (\alpha_2,0)$ or $(1,\beta_4)$ or $(-1,\beta_4)$,
  \item  $A^5_{\alpha_2}$ with $\alpha_2 \neq 0,1,-1$.
 \end{enumerate}
For anyone of these algebras $A$, there exists an unital $\K$ algebra $B_A=(V,\mu_{u,A})$ and linear endomorphisms $f_A,g_A$ such that the multiplication of $A$ is given by
$$\mu_A(X,Y)=\mu_{u,A}(f(X),g(Y)).$$
\end{proposition}
This unital algebra $B_A$ is called the unital heart of $A$. Since $B_A$ is unital, then (\cite{Pet}) it is an etale algebra, that is $B_A \otimes \widetilde{\K}=\widetilde{\K}^2$ where $\widetilde{\K}$ is the algebraic closure of $A$, or $B_A$ is isomorphic to the dual algebra defined by
$\mu_B(e_1,e_i)=\mu_B(e_i,e_1)=e_i, i=1,2$ and $\mu_B(e_2,e_2)=0$. To find this heart algebra we use the Kaplansky's Trick. If $A$ is regular, we consider $U$ and $V$ such that $L_U$ and $R_V$ are non singular and $f=L_U^{-1}$, $g=R_T^{-1}$. The multiplication $\mu_u$ of  the heart $B$ is $\mu_u(X,Y)=\mu(g(X),f(Y)$ and the identity of $B$ is $1_B=\mu(U,T)$.
\begin{enumerate}
  \item Let  be $A^1_{\alpha_2,\beta_2,\alpha_4,\beta_4}$. If $\alpha_2 \neq 1$ or $-1$ then $L_{e_1}$ and $R_{e_1}$ are not singular. In fact
  $$
   L_{e_1}=\begin{pmatrix}
   0   &  \alpha_2+1  \\
    1  &  \beta_2
\end{pmatrix}
, \ 
 R_{e_1}=\begin{pmatrix}
   0   &  \alpha_2-1  \\
    1  &  \beta_2
\end{pmatrix}
$$
Thus 
$$
  f=\ds \frac{-1}{\alpha_2 +1}\begin{pmatrix}
  \beta_2   &  -\alpha_2-1  \\
    -1  & 0
\end{pmatrix}
, \ 
g=   \ds \frac{-1}{\alpha_2 -1}\begin{pmatrix}
  \beta_2   &  -\alpha_2+1  \\
    -1  & 0
\end{pmatrix}
$$
Then the identity element of $B_A$ is $e_2$ and 
$$\mu_B(e_1,e_1)=\mu_A(g(e_1)g(e_1))=\frac{1}{\alpha_2^2-1}(\beta_2e_1-e_2)^2$$ and $B_A$ is etale.
If $\alpha_2=-1$, then we can take $U=e_2$ and $T=e_1$ as soon as $\alpha_4\beta_2 \neq 2 \beta_4$. If not we take $U=e_1+e_2$ and $T=e_1$. We have the same calcul for $\alpha_2=1$.
  \item Let be $A^2_{\alpha_1,\alpha_2,\alpha_4}$. This algebra is regular. If $\alpha_1 \neq 0$, then $L_{e_1}$ and $R_{e_1}$ are not singular
  and $B_A$ is etale.
\end{enumerate}

 \subsection{Case $\mu_a(e_1,e_2)=0$}

The multiplication $\mu$ is symmetric. The group of automorphisms of $\mu_a$ is $GL(2,\K)$.  Moreover the multiplication $\mu$ writes
$$
\left\{
\begin{array}{l}
\mu(e_1,e_1)=\alpha_1e_1+\beta _1e_2, \\
\mu(e_1,e_2)=\alpha_2e_1+\beta _2e_2, \\
\mu(e_2,e_1)=\alpha_2e_1+\beta _2e_2, \\
 \mu(e_2,e_2)=\alpha_4e_1+\beta _4e_2, \\
 \end{array}
 \right.
 $$

 $\bullet$ We assume that there exists two independent idempotent vectors. If $e_1$ and $e_2$ are these vectors, then
 $$\mu(e_1,e_1)=e_1, \ \mu(e_2,e_2)=e_2.$$ We obtain the following algebras
 $$
\left\{
\begin{array}{l}
\mu^{6}_{\alpha_2,\beta_2}(e_1,e_1)=e_1, \\
\mu^{6}_{\alpha_2,\beta_2}(e_1,e_2)=\alpha_2e_1+\beta _2e_2, \\
\mu^{6}_{\alpha_2,\beta_2}(e_2,e_2)=e_2. \\
 \end{array}
 \right.
 $$
 Remark that if any element is idempotent, thus $\mu(e_1,e_2)=\mu(e_2,e_1)=0.$ In fact
 $$\mu(e_1+e_2,e_1+e_2)=e_1+e_2=\mu(e_1,e_1)+\mu(e_2,e_2)+2\mu(e_1,e_2).$$
 In the general case, if $ae_1+be_2$ is an idempotent with $ab \neq 0$, then $a$ and $b$ satisfy the system
 $$\left\{
 \begin{array}{ll}
 a^2+2ab\alpha_2 = a\\
 b ^2+2ab\beta_2 =b.
 \end{array}
 \right.
 $$
 If $4\alpha_2\beta_2=1$, then the system has solutions as soon as $\alpha_2=\beta_2=\frac{1}{2}$. In this case we obtain the multiplication
 $\mu^{6}_{\frac{1}{2},\frac{1}{2}}$ and for any $a$, the vectors $ae_1+(1-a)e_2$ are idempotent. If $4\alpha_2\beta_2 \neq 1$, the vector
 $$v=\displaystyle\frac{1-2\alpha_2}{1-4\alpha_2\beta_2}e_1+\frac{1-2\beta_2}{1-4\alpha_2\beta_2}e_2
$$
is an idempotent and the only idempotents are $e_1, e_2$ and $v.$ The changes of basis
$\{ e_1,v\}$ or $\{ e_2, v \}$ do not simplify the number of independent parameters.

\medskip

 $\bullet$ We assume that there exists only one idempotent vector. If $e_1$ is this vector, thus $\mu(e_1,e_1)=e_1.$  If we consider a vector $v=xe_1+ye_2$ such that $\mu(v,v)=v,$ then $x$ and $y$ have to satisfy
  \begin{equation}
\label{toto}
\left\{
 \begin{array}{l}
     x^2+2xy\alpha_2+y ^2 \alpha_4=x,    \\
     2xy\beta_2+y ^2\beta_4=y.
\end{array}
\right.
\end{equation}
If we assume that $y \neq 0$, the second equation gives as soon as $\beta_2 \neq 0$,  $x=\ds \frac{1-y\beta_4}{2\beta_2}$ and thus
 \begin{equation}
 \label{equation1idem}
 y^2(\beta_4^2-4\alpha_2\beta_2\beta_4+4\beta_2^2\alpha_4)+y(4\alpha_2\beta_2+2\beta_2\beta_4-2\beta_4)+1-2\beta_2=0.
 \end{equation}
 Let us consider
 a change of basis which preserves $e_1,$ that is,
 \begin{equation}
 \left\{
 \begin{array}{l}
 \label{chtbase}
 e'_1=e_1, \\
  e'_2=be_1+de_2,
  \end{array}
  \right.
   \end{equation}
 with $d \neq 0.$ Since in this new basis we have $\beta'_4=2b\beta_2+d\beta_4$, we can find $b$ such that $\beta'_4=0$. Then we can assume that $\beta_4=0$.

  If  moreover $\alpha_2 \neq 0$, taking $d=\alpha_2^{-1}$, we obtain $\alpha'_2=1$ and we have the algebra
$$
\left\{
\begin{array}{l}
\mu(e_1,e_1)=e_1, \\
\mu(e_1,e_2)=e_1+\beta _2e_2, \\
 \mu(e_2,e_2)=\alpha_4e_1.\\
 \end{array}
 \right.
$$
 Equation (\ref{equation1idem})
simplifies as
\begin{equation}
\label{y}
y^2(4\beta_2^2\alpha_4)+4\beta_2y+1-2\beta_2=0.
\end{equation}
 If we assume that $\K$ is algebraically closed, then this equation has in general two roots. It has no root if $\beta_2=0$ which is excluded.  Then to have only one idempotent,  $0$ must be the only root  which is equivalent to $\alpha_4=0$ and $\beta_2=1/2$.  We obtain the following algebra
 $$
\left\{
\begin{array}{l}
\mu^{7}(e_1,e_1)=e_1, \\
\mu^{7}(e_1,e_2)=\ds e_1+ \frac{1}{2}e_2, \\
 \mu^{7}(e_2,e_2)=0.\\
 \end{array}
 \right.
 $$
 If $\K$ is not algebraically closed, then we have no idempotent other than $0$  if $\alpha_4=0$ and $\beta_2=1/2$ and we obtain the previous algebra $\mu^7$ or if $y^2(4\beta_2^2\alpha_4)+4\beta_2y+1-2\beta_2$ is irreducible in $\K$. We obtain
 $$
\left\{
\begin{array}{l}
\mu^{7}_R(e_1,e_1)=e_1, \\
\mu^{7}_R(e_1,e_2)=\ds e_1+ \beta_2e_2, \\
 \mu^{7}_R(e_2,e_2)=\alpha_4 e_1,\\
 \end{array}
 \right.
 $$
 with $y^2(4\beta_2^2\alpha_4)+4\beta_2y+1-2\beta_2$   irreducible in $\K$ (so $\alpha_4 \neq 0$).

 \noindent If $\alpha_2=0$ and if $\K$ is algebraically closed, we consider in the change of basis (\ref{chtbase}) defined above, $b=0$ and $d=\sqrt{\alpha_4}$ if $\alpha_4 \neq 0$:
 $$
\left\{
\begin{array}{l}
\mu(e_1,e_1)=e_1, \\
\mu(e_1,e_2)=\beta _2e_2, \\
 \mu(e_2,e_2)=e_1.\\
 \end{array}
 \right.
 $$
 There exits only one idempotent if and only if $\beta_2=1/2$. We obtain the following algebra
  $$
\left\{
\begin{array}{l}
\mu^{8}(e_1,e_1)=e_1, \\
\mu^{8}(e_1,e_2)=\ds \frac{1}{2}e_2, \\
 \mu^{8}(e_2,e_2)=e_1.\\
 \end{array}
 \right.
 $$
 If $\alpha_2=\alpha_4=0$,  we have only one idempotent if and only if $2\beta_2 \neq 1$. We obtain
  $$
\left\{
\begin{array}{l}
\mu^{9}(e_1,e_1)=e_1, \\
\mu^{9}(e_1,e_2)=\beta_2 e_2,  \ \ (\beta_2 \neq 1/2)\\
 \mu^{9}(e_2,e_2)=0.\\
 \end{array}
 \right.
 $$
Assume  $\K$  not algebraically closed and $\alpha_2=0$.  If the equation $d^2 \alpha_4$ has a root in $\K$, we find $\mu^8$. If not, let $\lambda_2 \in \K/(\K^*)^2$ such that $d^2\alpha_4= \lambda_2$. In this case we have only one idempotent if and only if  $(2\beta_2=1)$ or $(1-2\beta_2 \notin (\K^*) ^2)$. We obtain
 $$
\left\{
\begin{array}{l}
\mu_R^{8,1}(e_1,e_1)=e_1, \\
\mu_R^{8,1}(e_1,e_2)=\ds \frac{1}{2}e_2,\\
 \mu_R^{8,1}(e_2,e_2)=\lambda_2e_1,\\
 \end{array}
 \right.
 $$
 and
 $$
\left\{
\begin{array}{l}
\mu_R^{8,2}(e_1,e_1)=e_1, \\
\mu_R^{8,2}(e_1,e_2)=\beta_2 e_2, \ \ 1-2\beta_2 \notin  (\K^*) ^2,\\
 \mu_R^{8,2}(e_2,e_2)=\lambda_2 e_1.\\
 \end{array}
 \right.
 $$

 \medskip

 \noindent
 Assume now that $\beta_2=0$. Then (\ref{toto}) implies $y^2\beta_4=y.$
 If $\beta_4=0$, then $y=0$ and we have
 $$
\left\{
\begin{array}{l}
\mu(e_1,e_1)=e_1, \\
\mu(e_1,e_2)=\alpha_2 e_1,  \\
 \mu(e_2,e_2)=\alpha_4e_1.\\
 \end{array}
 \right.
 $$
 The change of basis $e'_1=e_1,e'_2=be_1+de_2$ gives $\alpha'_2=d\alpha_2,\alpha_4'=d^Ž\alpha_4$. We obtain
 $$
\left\{
\begin{array}{l}
\mu^{10}(e_1,e_1)=e_1, \\
\mu^{10}(e_1,e_2)= e_1,  \\
 \mu^{10}(e_2,e_2)=\alpha_4e_1.\\
 \end{array}
 \right.
 $$
 if $\alpha_2 \neq 0$.  Assume now that $\alpha_2=0$ and $\alpha_4 \neq 0$. If $\K$ is algebraically close, we obtain
 $$
\left\{
\begin{array}{l}
\mu^{11}(e_1,e_1)=e_1, \\
\mu^{11}(e_1,e_2)= 0,  \\
 \mu^{11}(e_2,e_2)=e_1,\\
 \end{array}
 \right.
 $$
$$
\left\{
\begin{array}{l}
\mu^{11}_R(e_1,e_1)=e_1, \\
\mu^{11}_R(e_1,e_2)= 0,  \\
 \mu^{11}_R(e_2,e_2)=\lambda_2e_1\\
 \end{array}
 \right.
 $$
 with $\lambda_2 \in \K/(\K^*)^2$. If $\alpha_4=0$,
 $$
\left\{
\begin{array}{l}
\mu^{12}(e_1,e_1)=e_1, \\
\mu^{12}(e_1,e_2)= 0,  \\
 \mu^{12}(e_2,e_2)=0\\
 \end{array}
 \right.
 $$
 $\bullet$ No vector is idempotent. If there exists $v$ with $\mu(v,v) \neq 0$, thus we can consider that $\mu(e_1,e_1)=e_2$ that is
 $$
\left\{
\begin{array}{l}
\mu(e_1,e_1)=e_2, \\
\mu(e_1,e_2)=\mu(e_2,e_1)=\alpha_2e_1+\beta _2e_2, \\
\mu(e_2,e_2)=\alpha_4e_1+\beta_4e_2.
 \end{array}
 \right.
 $$
 
1.  If $\alpha_4=0$, that is $\mu(e_2,e_2)=\beta_4e_2$, then the vector $e'_2=\beta_4^{-1}e_2$ is idempotent as soon as $\beta_4 \neq 0$. Then the hypothesis implies $\beta_4=0$.  Let be $v=xe_1+ye_2$. The equation $\mu(v,v)=v$ 
 is equivalent to:
 $$x^2e_2+2xy(\alpha_2e_1+\beta _2e_2)=2xy\alpha_2e_1+(x^2+2xy\beta _2)e_2=xe_1+ye_2.$$
that is 
 $$2xy\alpha_2=x,  \ x^2+2xy\beta _2=y.$$
 If $\alpha_2=0$, then $x=y=0$, and no elements are idempotent. We obtain the algebras, corresponding to $\beta_2 \neq 0$ or $\beta_2=0$
 $$
\left\{
\begin{array}{l}
\mu^{13}(e_1,e_1)=e_2, \\
\mu^{13}(e_1,e_2)=e_2, \\
\mu^{13}(e_2,e_2)=0.
 \end{array}
 \right.
 $$
 $$
\left\{
\begin{array}{l}
\mu^{14}(e_1,e_1)=e_2, \\
\mu^{14}(e_1,e_2)=0, \\
\mu^{14}(e_2,e_2)=0.
 \end{array}
 \right.
 $$
 If $\alpha_2 \neq 0$ and $y=(2\alpha_2)^{-1}$ then $x$ satisfies the equation
 \begin{equation}
\label{uu}
\displaystyle x^2+\left(\frac{\beta _2}{\alpha_2}\right)x-\frac{1}{2\alpha_2}=0.
\end{equation}
  If $\K$ is algebraically closed, such equation admits a non trivial solution. This is not compatible with our hypothesis.
Assume that  $\K$ is not algebraically closed. If $\beta_2 \neq 0$,  the change of basis $e'_1=\beta_2 ^{-1}e_1$ and $e'_2=\beta_2^{-2}e_2$ permits to consider $\beta_2=1$ and the (\ref{uu}) becomes
$$x^2+\frac{1}{\alpha_2}x-\frac{1}{2\alpha_2}=(x+\frac{1}{2\alpha_2})^2-\frac{1+2\alpha_2}{4\alpha_2^2}$$
 This equation has a non solution if $1+2\alpha_2 \notin (\K)^2$ where $(\K)^2=\{\lambda^2, \ \lambda \in \K\}.$
We obtain the algebras
     $$
\left\{
\begin{array}{l}
\mu^{14,1}_{R}(e_1,e_1)=e_2, \\
\mu^{14,1}_{R}(e_1,e_2)=\alpha_2e_1+e_2,   \ \ 2\alpha_2+1 \notin (\K)^2,\\
\mu^{14,1}_{R}(e_2,e_2)=0, \\
 \end{array}
 \right.
 $$
and
$$
\left\{
\begin{array}{l}
\mu^{14,2}_{R}(e_1,e_1)=e_2, \\
\mu^{14,2}_{R}(e_1,e_2)=\alpha_2e_1   \ \ 2\alpha_2 \notin (\K)^2,\\
\mu^{14,2}_{R}(e_2,e_2)=0. \\
 \end{array}
 \right.
 $$
 
 2. If  $\alpha_4 \neq 0$ the vector $v=xe_1+ye_2$ is idempotent if and only if
 $$ \left\{
 \begin{array}{l}
   2xy\alpha_2 +  y^2\alpha_4 =  x,    \\
    x^2 +2xy\beta_2+y^2\beta_4 = y.
\end{array}
 \right.
 $$
Then $x=\ds \frac{y^2\alpha_4}{1-2y\alpha_2}$. Let us note that   $1-2y\alpha_2 \neq 0$ because $1-2y\alpha_2=0$ implies $y^2\alpha_4=0$ that is $y=0$ and in this case $x=0$ and $v=0$. We deduce that $y$ is a root of the equation
$$
 (\ds \frac{y^2\alpha_4}{1-2y\alpha_2})^2 +2\ds \frac{y^2\alpha_4}{1-2y\alpha_2}y\beta_2+y^2\beta_4 - y=0$$
 that is
 $$-1+y(4\alpha_2+\beta_4)+y^2(2\alpha_4\beta_2-4\alpha_2^2-4\alpha_2\beta_4)+y^3(\alpha_4^2-4\alpha_2\alpha_4\beta_2+4\alpha_2^2\beta_4)=0.$$
If $\K$ is algebraically closed, this equation admits always a solution except if
$$ \left\{
 \begin{array}{l}
   4\alpha_2+\beta_4=0,   \\
   2\alpha_4\beta_2-4\alpha_2^2-4\alpha_2\beta_4=0,\\
   \alpha_4^2-4\alpha_2\alpha_4\beta_2+4\alpha_2^2\beta_4=0.
\end{array}
 \right.
 $$
 Then $\beta_4=-4\alpha_2, \alpha_4\beta_2=-6\alpha_2^2, \alpha_4^2=-8\alpha_2^3.$   We note that $\beta_2=0$ implies, if the characteristic of $\K$ is not $3$,  $\alpha_2=\alpha_4=0.$ From hypothesis, we can assume that $\beta_2 \neq 0$ and the change of basis $e'_1=ke_1,e'_2=k^2e_2$ which preserves the condition $e_1e_1=e_2$ changes $\beta_2$ in $k\beta_2$ and we can take $\beta_2=3$.  Then $\alpha_4=-2\alpha_2^2$,  $\alpha_4^2=4\alpha_2^4=-8\alpha_2^3$, then $\alpha_2=-2$ and $\alpha_4=4,\beta_4=8$ and we obtain
 the algebra
  $$
\left\{
\begin{array}{l}
\mu^{15}(e_1,e_1)=e_2, \\
\mu^{15}(e_1,e_2)=-2e_1+3e_2, \\
\mu^{15}(e_2,e_2)=-8e_1+8e_2.  \\
 \end{array}
 \right.
 $$
 Let us note that if the characteristic of $\K$ is $3$, then $\alpha_4\beta_2=0$ and $\beta_2=0$. This gives $\alpha_2(\alpha_2+\beta_4)=0$ and $\alpha_4^2+4\alpha_2^2\beta_4=0$. Since $\alpha_2=0$  implies $\alpha_4=0$ and $4\alpha_2+\beta_4=\alpha_2+\beta_4=0$ we obtain   $\beta_4=2\alpha_2$ and $ \alpha_4^2=2\alpha_2^2\beta_4=\alpha_2^3.$ By a change of basis we can take $\alpha_2=1$ and we obtain the algebra
 $$
\left\{
\begin{array}{l}
\mu^{15}_{(3)}(e_1,e_1)=e_2, \\
\mu^{15}_{(3)}(e_1,e_2)=e_1, \\
\mu^{15}_{(3)}(e_2,e_2)=e_1+2e_2.  \\
 \end{array}
 \right.
 $$ 
which correspond to $\mu_{15}$ in characteristic $3$.
 
If $\K$ is not algebraically closed, we have to consider all the algebras for which the polynomial
\begin{equation}
\label{ PA}
P_A(y)=-1+y(4\alpha_2+\beta_4)+y^2(2\alpha_4\beta_2-4\alpha_2^2-4\alpha_2\beta_4)+y^3(\alpha_4^2-4\alpha_2\alpha_4\beta_2+4\alpha_2^2\beta_4)
\end{equation}
has no root this is equivalent to say that $P_A$ is irreducible. If we consider the coefficient of $y^3$, that is $q_3(A)=\alpha_4^2-4\alpha_2\alpha_4\beta_2+4\alpha_2^2\beta_4$, it is equal to the discriminant of the determinant of the endomorphism$L_V$, that is $q_3(A)=Disc (\det (L_V))$. We deduce
\begin{proposition}
The algebra $A$ is regular if and only if $P_A(y)$ is strictly of degree $3$.
\end{proposition}

It remains to examine the case $\mu(v,v)=0$ for any $v$. That is
  $$
\left\{
\begin{array}{l}
\mu(e_1,e_1)=0, \\
\mu(e_1,e_2)=\alpha_2e_1+\beta_2e_2, \\
\mu(e_2,e_2)=0.  \\
 \end{array}
 \right.
 $$
 If $\alpha_2\beta_2 \neq 0$ we can find some idempotents. In all the others cases, we have no idempotent. We obtain
  $$
\left\{
\begin{array}{l}
\mu^{16}(e_1,e_1)=0, \\
\mu^{16}(e_1,e_2)=e_1, \\
\mu^{16}(e_2,e_2)=0,  \\
 \end{array}
 \right.
 $$
 and
  $$
\left\{
\begin{array}{l}
\mu^{17}(e_1,e_1)=0, \\
\mu^{17}(e_1,e_2)=0, \\
\mu^{17}(e_2,e_2)=0.  \\
 \end{array}
 \right.
 $$

 \begin{theorem}
 Any commutative $2$-dimensional algebra  over an algebraically closed field is isomorphic to one of the following

 $$\bullet \ \ \left\{
\begin{array}{l}
\mu^{6}(e_1,e_1)=e_1, \\
\mu^{6}(e_1,e_2)=\alpha_2e_1+\beta _2e_2, \\
\mu^{6}(e_2,e_2)=e_2. \\
 \end{array}
 \right. \ \
\left\{
\begin{array}{l}
\mu^{7}(e_1,e_1)=e_1, \\
\mu^{7}(e_1,e_2)=\ds e_1+ \frac{1}{2}e_2, \\
 \mu^{7}(e_2,e_2)=0.\\
 \end{array}
 \right. \ \  \left\{
\begin{array}{l}
\mu^{8}(e_1,e_1)=e_1, \\
\mu^{8}(e_1,e_2)=\ds \frac{1}{2}e_2, \\
 \mu^{8}(e_2,e_2)=e_1.\\
 \end{array}
 \right.
 $$
 $$ \bullet \   \left\{
\begin{array}{l}
\mu^{9}(e_1,e_1)=e_1, \\
\mu^{9}(e_1,e_2)=\beta_2 e_2,  \ \ (\beta_2 \neq 1/2),\\
 \mu^{9}(e_2,e_2)=0.\\
 \end{array}
 \right.
 \ \
\left\{
\begin{array}{l}
\mu^{10}(e_1,e_1)=e_1, \\
\mu^{10}(e_1,e_2)= e_1,  \\
 \mu^{10}(e_2,e_2)=\alpha_4e_1.\\
 \end{array}
 \right.
 \ \
\left\{
\begin{array}{l}
\mu^{11}(e_1,e_1)=e_1, \\
\mu^{11}(e_1,e_2)= 0,  \\
 \mu^{11}(e_2,e_2)=e_1.\\
 \end{array}
 \right.
$$
 $$ \bullet \
\left\{
\begin{array}{l}
\mu^{12}(e_1,e_1)=e_1, \\
\mu^{12}(e_1,e_2)= 0,  \\
 \mu^{12}(e_2,e_2)=0.\\
 \end{array}
 \right.
\ \
Ê\left\{
\begin{array}{l}
\mu^{13}(e_1,e_1)=e_2, \\
\mu^{13}(e_1,e_2)=e_2, \\
\mu^{13}(e_2,e_2)=0.
 \end{array}
 \right.
 \ \
\left\{
\begin{array}{l}
\mu^{14}(e_1,e_1)=e_2, \\
\mu^{14}(e_1,e_2)=0, \\
\mu^{14}(e_2,e_2)=0.
 \end{array}
 \right.
 $$
 $$ \bullet  
\left\{
\begin{array}{l}
\mu^{15}(e_1,e_1)=e_2, \\
\mu^{15}(e_1,e_2)=-2e_1+3e_2, \\
\mu^{15}(e_2,e_2)=-8e_1+8e_2.  \\
 \end{array}
 \right.
\ \ \left\{
\begin{array}{l}
\mu^{16}(e_1,e_1)=0, \\
\mu^{16}(e_1,e_2)=e_1, \\
\mu^{16}(e_2,e_2)=0.  \\
 \end{array}
 \right.
 \ \
\left\{
\begin{array}{l}
\mu^{17}(e_1,e_1)=0, \\
\mu^{17}(e_1,e_2)=0, \\
\mu^{17}(e_2,e_2)=0.  \\
 \end{array}
 \right.
 $$
 If $\K$ is not algebraically closed, we have also the following algebras where $\lambda_2 \in \K/(\K^*)^2$:
 $$ \bullet \ \
\left\{
\begin{array}{l}
\mu_R^{8,1}(e_1,e_1)=e_1, \\
\mu_R^{8,1}(e_1,e_2)=\ds \frac{1}{2}e_2,\\
 \mu_R^{8,1}(e_2,e_2)=\lambda_2e_1.\\
 \end{array}
 \right. \ \
 \left\{
\begin{array}{l}
\mu_R^{8,2}(e_1,e_1)=e_1, \\
\mu_R^{8,2}(e_1,e_2)=\beta_2 e_2, \ 1-2\beta_2 \notin  (\K^*) ^2,\\
 \mu_R^{8,2}(e_2,e_2)=\lambda_2 e_1.\\
 \end{array}
 \right.
\ \
\left\{
\begin{array}{l}
\mu^{11}_R(e_1,e_1)=e_1, \\
\mu^{11}_R(e_1,e_2)= 0,  \\
 \mu^{11}_R(e_2,e_2)=\lambda_2e_1. \\
 \end{array}
 \right.
 $$

   $$ \bullet
\left\{
\begin{array}{l}
\mu^{14,1}_{R}(e_1,e_1)=e_2, \\
\mu^{14,1}_{R}(e_1,e_2)=\alpha_2e_1+e_2,   \ \ 2\alpha_2+1 \notin \K^2,\\
\mu^{14,1}_{R}(e_2,e_2)=0. \\
 \end{array}
 \right.
 \ \
\left\{
\begin{array}{l}
\mu^{14,2}_{R}(e_1,e_1)=e_2, \\
\mu^{14,2}_{R}(e_1,e_2)=\alpha_2e_1   \ \ 2\alpha_2+1 \notin \K^2,\\
\mu^{14,2}_{R}(e_2,e_2)=0. \\
 \end{array}
 \right.
 $$
 $$ \bullet
\left\{
\begin{array}{l}
\mu^{15,1}_{R}(e_1,e_1)=e_2, \\
\mu^{15,1}_{R}(e_1,e_2)=\alpha_2e_1+\beta_2e_2,   P_A(y) \ {\rm without \ roots} \\
\mu^{15,1}_{R}(e_2,e_2)=\alpha_4e_1+\beta_4e_2. \\
 \end{array}
 \right.
$$

 \end{theorem}

Let us examine the property of regularity for these algebras. Since they are commutative, the left and right regularity are equivalent notions. Computing directly the determinant of the operator $L_{xe_1+ye_2}$ we deduce in the case $\K$ algebraically closed:
\begin{enumerate}
  \item The algebras $A^6=(V,\mu_6)$,  $A^7=(V,\mu_7)$, $A^8=(V,\mu_8)$, $A^{10}=(V,\mu_{10})$,  $A^{15}=(V,\mu_{15})$ are regular,
  \item $A^9=(V,\mu_9)$ is regular if $\beta_2 \neq 0$, 
  \item The algebras $A^{11}=(V,\mu_{11})$, $A^{12}=(V,\mu_{12})$, $A^{13}=(V,\mu_{13})$, $A^{14}=(V,\mu_{14})$, $A^{16}=(V,\mu_{16})$ and $A^{17}=(V,\mu_{17})$ are bisingular.
\end{enumerate}

 \section{Algebras over a field of characteristic $2$}
Let $\F$ be a field of characteristic $2$. Assume that $\F=\F_2$. If $A$ is a $2$-dimensional $\F$-algebra and if $\{e_1,e_2\}$ is a basis of $A$, then the values of the different products
belong to $\{e_1,e_2,e_1+e_2\}$. If $f$ is an isomorphism of $A$, it is represented in the basis $\{e_1,e_2\}$ by one of the following matrices
$$
M_1=\left(
\begin{array}{ll}
1&0\\
0&1
\end{array}
\right)
,\ M_2=\left(
\begin{array}{ll}
0&1\\
1&0
\end{array}
\right)
,\ M_3=\left(
\begin{array}{ll}
1&0\\
1&1
\end{array}
\right)
$$
$$
M_4=\left(
\begin{array}{ll}
1&1\\
0&1
\end{array}
\right)
,\ M_5=\left(
\begin{array}{ll}
0&1\\
1&1
\end{array}
\right)
,\ M_6=\left(
\begin{array}{ll}
1&1\\
1&0
\end{array}
\right)
$$
Each  of these matrices corresponds to a permutation of the finite set
$\{e_1,e_2,e_3=e_1+e_2\}$. If fact we have the correspondance:
$$
\begin{array}{|c|c|}
\hline
GL(A) & \Sigma_3 \\
\hline
M_1 & Id \\
\hline
M_2 & \tau_{12} \\
\hline
M_3 & \tau_{13} \\
\hline
M_4 & \tau_{23} \\
\hline
M_5 & c \\
\hline
M_6 & c^2 \\
\hline
\end{array}
$$
where $\tau _{ij}$ is the transposition between $i$ and $j$ and $c$ the cycle $(231)$. In fact, the matrix  $M_2$
corresponds to the linear  transformation
$f_2(e_1)=e_2, \ f(e_2)=e_1$ and in the set $(e_1,e_2,e_3)$ we have the transformation whose image is $(e_2,e_1,e_3)$ that is
the transposition $\tau _{12}$. The matrix $M_3$ corresponds to the linear transformation $f_2(e_1)=e_1+e_2, \ f(e_2)=e_2$ which corresponds to the permutation
$(e_3,e_2,e_1)$ that is $\tau _{13}$. For all other matrices we have similar results. We deduce
\begin{theorem}
There is a one-to-one correspondance between the change of  $\F$-basis in $A$ and the group $\Sigma_3$.
\end{theorem}
If we want to classify all these products of $A$, we have to consider all the possible results of these products and to determine the orbits of
the action of $\Sigma_3$.  More precisely the product $\mu(e_i,e_j)$ is in values in the set $(e_1,e_2,e_3=e_1+e_2)$. If we write
$\mu(e_i,e_j)=ae_1+be_2+ce_3$, thus the matrix $(a,b,c)$ is one of the following
$$R_0=(0,0,0)=0,R_1=(1,0,0), \ R_2=(0,1,0), \ R_3=(0,0,1).$$
Let us consider the following sequence
$$(\mu(e_1,e_1),\mu(e_1,e_2),\mu(e_2,e_1),\mu(e_2,e_2), \mu(e_1,e_3),\mu(e_2,e_3),\mu(e_3,e_1),\mu(e_3,e_2),\mu(e_3,e_3)).$$
As $\mu(e_1,e_3)=\mu(e_1,e_1+e_2)$, if $\mu(e_1,e_1)=R_i$ and $\mu(e_1,e_2)=R_j$ then $\mu(e_1,e_3)=R_i+R_j$ with the relations
$$R_i+R_i=0, \ R_i+R_j=R_k,$$
for $i,j,k$ all different and non zero.
Thus the four  first terms of this sequence determine all the other terms. More precisely, such a sequence writes
$$(R_i,R_j,R_k,R_l,R_i+R_j,R_k+R_l,R_i+R_k,R_j+R_l,R_i+R_j+R_k+R_l).$$

\noindent{\bf Consequence.} We have $4^4=256$ sequences, each  of these sequences corresponds to a $2$-dimensional $\F$-algebra.

\medskip

Let us denote by $S$ the set of these sequences. We have an action of $\Sigma_3$ on $S$: if $\sigma \in \Sigma_3$ and $s \in S$, thus
$s'=\sigma s$ is the sequence
$$
\begin{array}{l}
(\mu(e_{\sigma (1)},e_{\sigma (1)}),\mu(e_{\sigma (1)},e_{\sigma (2)}),\mu(e_{\sigma (2)},e_{\sigma (1)}),
\mu(e_{\sigma (2)},e_{\sigma (2)}), \mu(e_{\sigma (1)},e_{\sigma (3)}),\mu(e_{\sigma (2)},e_{\sigma (3)}),\\
\mu(e_{\sigma (3)},e_{\sigma (1)}),\mu(e_{\sigma (3)},e_{\sigma (2)}),\mu(e_{\sigma (3)},e_{\sigma (3)}))
\end{array}$$
with $\mu(e_{\sigma (i)},e_{\sigma (j)})=R_{\sigma ^{-1}(k)}$ when $\mu(e_i,e_j)=R_k$ and $R_k \neq 0.$ If $R_k=0$, then $\mu(e_{\sigma (i)},e_{\sigma (j)})=0.$
The classification of the $2$-dimensional $\F$-algebras corresponds to the determination of the orbits of this action. Recall that the subgroups of $\Sigma_3$ are $G_1=\{Id\},G_2= \{Id,\tau_{12}\},G_3=\{Id,\tau_{13}\},
G_4=\{Id,\tau_{23}\}, G_5=\{Id,c,c^2\},G_6=\Sigma_3$.
\begin{enumerate}
\item The isotropy subgroup is $\Sigma_3$. In this case we have the following sequence (we write only the $4$ first terms which determine the algebras:
    $$
    \begin{array}{l}
    s_1=(0,0,0,0)\\
    s_2=(R_1,R_3,R_3,R_2)
    \end{array}
    $$
    Recall that $\mu(e_1,e_1)=R_1$ means $\mu(e_1,e_1)=e_1$, $\mu(e_1,e_2)=R_3$ means  $\mu(e_1,e_2)=e_3$
and so on.
\item The isotropy subgroup is $G_5=\{Id,c,c^2\}$
We have only one orbit
$$
\begin{array}{|l|l|}
\hline
s & \mathcal{O}(s)\\
\hline
s_3=(R_3,R_2,R_2,R_1) & s_3, (R_2,R_1,R_1,R_3)\\
\hline
\end{array}
$$
\item The isotropy subgroup is of order $2$.
$$
\begin{array}{|l|l|}
\hline
s & \mathcal{O}(s)\\
\hline
s_4=(0,R_1,R_2,0) & s_4, (R_1,R_3,R_2,0), (0,R_1,R_3,R_2)\\
\hline
s_5=(0,R_2,R_1,0) & s_5,(R_1,R_2,R_3,0),(0,R_3,R_1,R_2) \\
\hline
s_6 = (0,R_3,R_3,0) & s_6, (0,R_1,R_1,0),(0,R_2,R_2),0)\\
\hline
s_7=(R_1,0,0,R_2) & s_7,(R_1,R_2,R_2,R_2),(R_1,R_1,R_1,R_2)\\
\hline
s_8=(R_1,R_1,R_2,R_2) & s_8,(0,R_1,0,R_2),(R_1,0,R_2,0)\\
\hline
s_9=(R_1,R_2,R_1,R_2) & s_9,(0,0,R_1,R_2),(R_1,R_2,0,0,)\\
\hline
s_{10}=(R_2,0,0,R_1) & s_{10},(R_1,R_3,R_3,R_3),(R_3,R_3,R_3,R_1)\\
\hline
s_{11}=(R_2,R_1,R_2,R_1) & s_{11},(0,0,R_1,R_3),(R_3,R_2,0,0)\\
\hline
s_{12}=(R_2,R_2,R_1,R_1) & s_{12},(0,R_1,0,R_3),(R_3,0,R_2,0)\\
\hline

s_{13}=(R_2,R_3,R_3,R_1)& s_{13},(R_1,R_2,R_2,R_3),(R_3,R_1,R_1,R_2)\\
\hline
s_{14}=(R_3,0,0,R_3) & s_{14},(0,R_1,R_1,R_1),(R_2,R_2,R_2,0)\\
\hline
s_{15}=(R_3,R_1,R_2,R_3) & s_{15},(R_1,R_2,R_3,R_1),(R_2,R_3,R_1,R_3)\\
\hline
s_{16}=(R_3,R_2,R_1,R_3) & s_{16},(R_1,R_3,R_2,R_1),(R_2,R_1,R_3,R_2)\\
\hline
s_{17}=(R_3,R_3,R_3,R_3) & s_{17}, (0,0,0,R_1),(R_2,0,0,0)\\
\hline
\end{array}
$$
\item The isotropy subgroup is trivial. In this case any orbit contains $6$ elements. As there are $256-46=210$ elements having $\Sigma_3$ as isotropy group, we deduce that we have $35$ distinguished non isomorphic classes.
\end{enumerate}
{\bf Conclusion} We have $52$ classes of non isomorphic algebras of dimension $2$ on the field $F_2$.

\section{Applications : $2$-dimensional $G$-associative and Jordan algebras}

 \subsection{$G$-associative commutative algebras}

 The notion of $G$-associativity has been defined in \cite{GozeRemm07}. Let $G$ be a subgroup of the symmetric group $\Sigma_3$. An algebra whose multiplication is denoted by $\mu$ is $G$-associative if we have
 $$\sum_{\sigma \in G}\varepsilon(\sigma)\mu(\mu(x_{\sigma(i)},x_{\sigma(j)}),x_{\sigma(k)})-\mu(x_{\sigma(i)},\mu(x_{\sigma(j)},x_{\sigma(k)})=0.
 $$
 where $\varepsilon(\sigma)$ is the signum of the permutation $\sigma$. Since we assume that $\mu$ is commutative, all these notions are trivial or coincide with the simple associativity. Now, if the algebra is of dimension $2$, then the associativity is completely determined by the identities
 $$\mu(\mu(e_1,e_1),e_2)-\mu(e_1,\mu(e_1,e_2)=0, \ \ \mu(\mu(e_1,e_2),e_2)-\mu(e_1,\mu(e_2,e_2)=0.$$
 We deduce that the only associative commutative $2$-dimensional algebras are
 \begin{itemize}
  \item $\mu^6$ for $(\alpha_2,\beta_2) \in \{(0,1),(1,0),(0,0)\},$
  \item $\mu^9$ for $\beta_2=0$ or $1$,
  \item $\mu^{12}, \mu^{16}, \mu^{17}.$
  \item if $\K=\R$: $\mu^8_R$ for $\beta_2=1$ and $\lambda=-1$.
\end{itemize}
 We find again the classical list (see for example \cite{GozeRemm03}).

  \subsection{$G$-associative noncommutative algebras}
Let us consider now the noncommutative case. From Theorem \ref{nc}, the multiplication $\mu$ is isomorphic to some
$\mu^i, i=1,\cdots,5$ (we consider here that $\K$ is algabraically closed). Let $A_\mu$ be the associator of $\mu$, that is $A_\mu=\mu \circ (\mu \otimes Id)- \mu \circ (Id \otimes \mu)$ and $\mu$ is associative if and only if $A_\mu =0$. The examination of this list  allows to find the classification of the $2$-dimensional noncommutative associative algebras: these algebras are isomorphic to one of the following
\begin{enumerate}
  \item $\mu^4_{-1,-2}$ that is $\left\{ \begin{array}{ l}
  e_1e_1=0,   \\
  e_1e_2=0,\\
  e_2e_1=-2e_1,\\
     e_2e_2=-2e_2.
\end{array}
\right.
$
  \item $\mu^4_{1,2}$ that is $\left\{ \begin{array}{ l}
  e_1e_1=0,   \\
  e_1e_2=2e_1,\\
  e_2e_1=0,\\
     e_2e_2=2e_2.
\end{array}
\right.
$
\end{enumerate}

\medskip

Now, for any nonassociative algebra, we examine the $G_i$-associativity. Note that all these algebras are Lie-admissible, that is $\Sigma_3$-associative. We focuse essentially on the $G_2$-associativity, $G_2=\{Id,\tau_{12}\}$, because we deduce immediately the affine structures on the associated Lie algebra $\mu_a$. Then we compute for any algebra $A_\mu(e_1,e_2,e_1)-A_\mu(e_2,e_1,e_1)$ and $A_\mu(e_1,e_2,e_2)-A_\mu(e_2,e_1,e_2)$. We deduce that $\mu^1_{\alpha_2,\beta_2,\alpha_4,\beta_4}$ is $G_2$-associative if and only if $\beta_2=\alpha_4=0$ and $\alpha_2=-1,\beta_4=-4$. The algebras $\mu^2$ and $\mu^3$ are never $G_2$-associative, $\mu^4_{\alpha_2,\beta_4}$ is $G_2$-associative for $\alpha_2=-1$ or $(\beta_4=\alpha_2-1)$. Likewise, $\mu^5_{\alpha_2}$ is $G_2$-associative for $\alpha_2=-1$ or $\alpha_2=1.$
\begin{proposition}
Any $2$-dimensional noncommutative $G_2$-associative algebra is isomorphic to one of the following
\begin{enumerate}
  \item $\mu^4_{-1,-2}$ or  $\mu^4_{1,2}$, that is $\mu$ is associative,
   \item $\mu^1_{-1,0,0,-4}$ that is $\left\{ \begin{array}{ l}
     e_1e_1=e_2,   \\
  e_1e_2=0,\\
  e_2e_1=-2e_1,\\
     e_2e_2=-4e_2.
\end{array}
\right.
$
  \item $\mu^4_{-1,\beta_4}$ that is $\left\{ \begin{array}{ l}
  e_1e_1=0,   \\
  e_1e_2=0,\\
  e_2e_1=-2e_1,\\
     e_2e_2=\beta_4e_2.
\end{array}
\right.
$
 \item $\mu^4_{\alpha_2,\alpha_2+1}$ that is $\left\{ \begin{array}{ l}
  e_1e_1=0,   \\
  e_1e_2=(\alpha_2+1)e_1,\\
  e_2e_1=(\alpha_2-1)e_1,\\
     e_2e_2=(\alpha_2+1)e_2.
\end{array}
\right.
$
 \item $\mu^5_{1}$ that is $\left\{ \begin{array}{ l}
  e_1e_1=0,   \\
  e_1e_2=2e_1,\\
  e_2e_1=0,\\
     e_2e_2=e_1+2e_2.
\end{array}
\right.
$
 \item $\mu^5_{-1}$ that is $\left\{ \begin{array}{ l}
  e_1e_1=0,   \\
  e_1e_2=0,\\
  e_2e_1=-2e_1,\\
     e_2e_2=e_1-2e_2.
\end{array}
\right.
$
\end{enumerate}
\end{proposition}

 \subsection{Jordan algebras}
In a Jordan algebra, the multiplication $\mu$ satisfies:
       $$
\left\{
\begin{array}{l}
\mu(v,w)=\mu(w,v)\\
\mu(\mu(v,w),\mu(v,v)) = \mu(v,\mu(w,\mu(v,v))
 \end{array}
 \right.
 $$
 for all $v,w$.
 We assume in this section that $\K$ is algebraically closed and that the Jordan algebra are of dimension $2$. Thus the multiplication $\mu$ is isomorphic to $\mu_i$ for $i=11, \cdots, 16.$
To simplify the notation, we will write $vw$ in place of $\mu(v,w)$. If $v$ is an idempotent, thus $v^2=v$ and the Jordan identity gives
$$v(vw)=v(vw)$$
for any $w$, that is, this identity is always satisfied.
\begin{lemma}
If $v_1$ and $v_2$ are idempotent vectors, thus
$$(v_1v_2)((v_1+v_2)w)=(v_1+v_2)((v_1v_2)w)$$
for any $w$.
\end{lemma}
\pf In the Jordan identity, we replace $v$ by $v_1+v_2$. We obtain
$$v_1^2(v_2w)+2(v_1v_2)((v_1+v_2)w)+v_2^2(v_1w)=v_1(v_2^2w)+v_2(v_1^2w)+2(v_1+v_2)((v_1v_2)w).$$
Since $v_1$ and $(v_2)$ are idempotent, this equation reduces
$$(v_1v_2)((v_1+v_2)w)=(v_1+v_2)((v_1v_2)w).$$
 \medskip
 \begin{proposition}
 If $v_1$ and $v_2$ are idempotent vectors such that $v_1v_2$ and $v_1+v_2$ are independent, thus the Jordan algebra is associative.
 \end{proposition}
 \pf Let $x$ and $y$ be two vectors of the algebra. Thus, by hypothesis, $x=x_1v_1v_2+x_2(v_1+v_2)$ and $y=y_1v_1v_2+y_2(v_1+v_2)$. Thus
 $$x(yw)=x_1y_1(v_1v_2)((v_1v_2)w)+(x_1y_2+x_2y_1)(v_1v_2)((v_1+v_2)w)+x_2y_2(v_1+v_2)((v_1+v_2)w).$$
 and
 $$x(yw)=y(xw).$$
 By commutativity we obtain
 $$x(yw)=x(wy)=y(xw)=(xw)y$$
 this proves that the algebra is associative.

 $\bullet$
 If $\mu$ is given by
  $$
\left\{
\begin{array}{l}
\mu(e_1,e_1)=e_1, \\
\mu(e_1,e_2)=\alpha_2e_1+\beta _2e_2, \\
\mu(e_2,e_2)=e_2 \\
 \end{array}
 \right.
 $$
 the Jordan algebra admits two idempotents $e_1$ and $e_2$. Since $e_1e_2=\alpha_2e_1+\beta _2e_2$, the vectors $e_1e_2$ and $e_1+e_2$ are independent if and only if $\alpha_2 \neq \beta_2$. In this case the algebra can be associative and we obtain the following associative Jordan algebra corresponding to
 \begin{enumerate}
 \item $\alpha_2=1$, $\beta_2=0$
 \item  $\alpha_2=0$, $\beta_2=1$
 \end{enumerate}
 These Jordan algebras are isomorphic. This gives the following Jordan algebra
   $$
   J_1=
\left\{
\begin{array}{l}
e_1e_1=e_1, \\
e_1e_2=e_2e_1=e_2\\
e_2e_2=e_2. \\
 \end{array}
 \right.
 $$
  If $e_1e_2$ and $e_1+e_2$ are dependent, that is $e_1e_2=\lambda(e_1+e_2)$, then $\lambda=-1$ or $\frac{1}{2}$ or $0$.  If  $e_1e_2=0$, the product is not a Jordan product.  If  $\lambda=-1$ the product is never a Jordan product. If  $\lambda=\frac{1}{2}$, we obtain the following Jordan algebra
 $$
   J_2=
\left\{
\begin{array}{l}
e_1e_1=e_1, \\
e_1e_2=e_2e_1=\frac{1}{2}(e_1+e_2)\\
e_2e_2=e_2. \\
 \end{array}
 \right.
 $$

 \medskip

 \medskip

 $\bullet$   $\mu$ is given by
 $$
\left\{
\begin{array}{l}
\mu(e_1,e_1)=e_1, \\
\mu(e_1,e_2)=\beta_2e_2, \\
\mu(e_2,e_2)=0. \\
 \end{array}
 \right.
 $$
This product is  a Jordan product if $\beta_2=1$ or $0$. We obtain
   $$
   J_3=
\left\{
\begin{array}{l}
e_1e_1=e_1, \\
e_1e_2=e_2e_1=e_2\\
e_2e_2=0. \\
 \end{array}
 \right.
 , \ \
    J_4=
\left\{
\begin{array}{l}
e_1e_1=e_1, \\
e_1e_2=e_2e_1=0\\
e_2e_2=0. \\
 \end{array}
 \right.
 $$
 $\bullet$  If $\mu=\mu_{11}$
 we have also a Jordan structure  $$
      J_5=
\left\{
\begin{array}{l}
e_1e_1=e_2 \\
e_1e_2=e_2e_1=0 \\
e_2e_2=0 \\
 \end{array}
 \right.
 $$

  \medskip

 $\bullet$ $\mu=0$, we have the trivial Jordan algebra.

 \medskip

 $\bullet$ If $\K$ is not algebraically closed, we consider
 $$
  \left\{
\begin{array}{l}
\mu_R^{8,2}(e_1,e_1)=e_1, \\
\mu_R^{8,2}(e_1,e_2)=\beta_2 e_2, \ \ 1-2\beta_2 \notin  (\K^*) ^2,\\
 \mu_R^{8,2}(e_2,e_2)=\lambda e_1,\\
 \end{array}
 \right.
$$
 We obtain a Jordan structure
 $$
      J_6=
\left\{
\begin{array}{l}
e_1e_1=e_1 \\
e_1e_2=e_2e_1=e_2 \\
e_2e_2=\lambda e_1. \\
 \end{array}
 \right.
 $$
We find the list established in \cite{Ancochea}.
 
 \subsection{$2$-dimensional Hom-algebra}
The notion of Hom-algebra was introduced to generalize those of Hom-Lie algebra which appeared naturally when we are interested by the notion of $q$-derivation on the Witt algebra. In dimension $2$, this notion is equivalent to the classical notion of Lie algebra. In dimension $3$, we have shown that any skew-symmetric algebra is a Hom-Lie algebra. Then our interest concerns Hom-associative algebra, that is algebra $A=(V,\mu)$ such that there exists $f \in End(V)$ satisfying the Hom-Ass identity:
$$\mu(\mu(X,Y),f(Z))=\mu(f(X),\mu(Y,Z))$$
for any $X,Y,Z \in V$. Using previous notations, we consider the algebras $A^{(Id,f)}$ and its opposite $A^{(f,Id)}$. Their multiplication law are respectively defined by
$$\mu_{R,f}(X,Y)=\mu(X,f(Y)), \ \ \mu_{L,f}(X,Y)=\mu(f(X),Y)$$
and the Hom-Ass identity can be written:
$$\mu_{R,f} \circ (\mu \otimes Id)- \mu_{L,f} \circ (Id \circ \mu)=0.$$
Assume now that the algebra $A$ is regular. In this case, assuming that the field is algebraically closed, there exists an unital algebra whose product is denoted  $X \cdot Y$ and two endomorphisms  $u$ and $v$  of $V$ such that
$$\mu(X,Y)=u(X)\cdot v(Y).$$
Then
$$\mu_{R,f}(X,Y)= u(X) \cdot v \circ f (Y), \ \ \mu_{L,f}(X,Y)=u \circ f (X) \cdot v(Y).$$
Then the Hom-Ass identity becomes
$$ u (u(X) \cdot v(Y)) \cdot v \circ f (Z)- u \circ f (X) \cdot (v (u(Y) \cdot v(Z))=0.$$
Maybe, it is better to look the Hom-Ass identity from the previous list. Assume that $A$ is non commutative.
\begin{enumerate}
\item $A=A^1_{\alpha_2,\beta_2,\alpha_4,\beta_4}=(V,\mu^1).$
 Let $f$ be an endomorphism of $V$ satisfying the Hom-Ass identity. To simplify notations we write $XY$ for $\mu(X,Y)$ and $[X,Y]$ for $\mu_a(X,Y)$.  We have in particular
 $$(e_1e_1)f(e_1)-f(e_1)(e_1e_1))=[e_2,f(e_1)]=0.$$
 We deduce $f(e_1)=a e_2.$ Likewise we have  $[e_2e_2,f(e_2)]=0$ and $f(e_2)=k(\alpha_4e_1+\beta_4e_2).$ Other identities give : \begin{enumerate}
  \item $(e_1e_2)f(e_1)-f(e_1)(e_2e_1)=0$ implies $a=0$ or $e_2e_2=0$.
  \item If $a=0$, then $(e_2e_1)f(e_1)-f(e_2)(e_1e_1)=0$ implies $f(e_2)e_2=0$ and $(e_1e_1)f(e_2)-f(e_1)(e_1e_2)=0$ implies $e_2f(e_2)=0$. Then  $[e_2,f(e_2)]=0$ and $f(e_2)=ke_2$.  This gives $0=f(e_2)e_2=be_2e_2$ that is $f=0$ or $e_2e_2=0$. But we have seen that $f(e_2)=k(\alpha_4e_1+\beta_4e_2)$, then in all the cases, $f=0$. 
  \item If $a \neq 0$, then $e_2e_2=0$ and $f(e_2)=0$. We deduce that $(e_1e_2)f(e_1)-f(e_1)(e_2e_1)=0$  implies $\alpha_2=\beta_2=0$. Thus $(e_2e_1)f(e_1)-f(e_2)(e_1e_1)=-a(e_1e_2)=-ae_1=0$ and $a=0$.
\end{enumerate}
We deduce that the algebra $A_{\alpha_2,\beta_2,\alpha_4,\beta_4}$ is not a Hom-associative algebra. 
\item $A=A^2_{\alpha_1,\alpha_2,\alpha_4}.$ With similar simple computation we can look that also this algebra is not a Hom-Ass algebra.
\item $A=A^3_{\alpha_4,\beta_4}.$ In this case also, if we compute $(e_1e_1)f(e_1)-f(e_1)(e_1e_1)=[e_1,f(e_1)]=0$, we obtain $f(e_1)=k_1e_1$. Also we have $(e_1e_2)f(e_1)-f(e_1)(e_2e_1)=2k_1e_1=0$ and $f(e_1)=0.$ We deduce $e_1f(e_2)=0$ and $f(e_2)e_1=0$ and $f(e_2)=0.$ Thus $f=0$ and $A^3$ is not a Hom-associative algebra.
\item  $A=A^4_{\alpha_2,\beta_4}.$ If $\beta_4 \neq 0$, then the Hom-Ass condition implies $\alpha_2=1$ or $-1$. We obtain the following Hom-Ass algebras:
$$
\left\{
\begin{array}{l}
\mu^4_{1,\beta_4}(e_1,e_1)=0, \\
\mu^4_{1,\beta_4}(e_1,e_2)=2e_1, \\
\mu^4_{1,\beta_4}(e_2,e_1)=0, \\
 \mu^4_{1,\beta_4}(e_2,e_2)=\beta_4 e_2, \\
 \end{array}
 \right.
 , \ \  \left\{
\begin{array}{l}
\mu^4_{-1,\beta_4}(e_1,e_1)=0, \\
\mu^4_{-1,\beta_4}(e_1,e_2)=0, \\
\mu^4_{-1,\beta_4}(e_2,e_1)=-2e_1, \\
 \mu^4_{-1,\beta_4}(e_2,e_2)=\beta_4 e_2. \\
 \end{array}
 \right.
$$
In each of these two cases, $f$ is a diagonal endomorphism.  These algebras are for $\beta_4 \neq 2$ or $-2$, not associative.
\item $A=A^5_{\alpha_2}$. If $\alpha_2=0$, any linear endomorphism with values in $\K\{e_1\}$ satisfies the Hom-Ass identity. Then the following algebra is Hom-associative
$$
\left\{
\begin{array}{l}
\mu^5_{0}(e_1,e_1)=0, \\
\mu^5_0(e_1,e_2)=e_1, \\
\mu^5_{0}(e_2,e_1)=-e_1, \\
 \mu^5_{0}(e_2,e_2)= e_1. \\
 \end{array}
 \right.
$$
Assume now that $\alpha_2 \neq 0$. If $\alpha_2 \neq \pm 1$, then any endomorphism satisfying the Hom-Ass identity is trivial.
If $\alpha_2=1$ or $-1$, we have non trivial solution and the following algebras are Hom-associative algebras
$$
\left\{
\begin{array}{l}
\mu^5_{-1}(e_1,e_1)=0, \\
\mu^5_{-1}(e_1,e_2)=0, \\
\mu^5_{-1}(e_2,e_1)=-2e_1, \\
 \mu^5_{-1}(e_2,e_2)= e_1-2e_2. \\
 \end{array}
 \right. , \ \  \left\{
\begin{array}{l}
\mu^5_{1}(e_1,e_1)=0, \\
\mu^5_{1}(e_1,e_2)=2e_1, \\
\mu^5_{1}(e_2,e_1)=0, \\
 \mu^5_{1}(e_2,e_2)= e_1+2e_2. \\
 \end{array}
 \right. 
$$
with $f=\begin{pmatrix}
    -4x  &  x  \\
     0 & -2x 
\end{pmatrix}$
in the first case and $f=\begin{pmatrix}
    4x  &  x  \\
     0 & 2x 
\end{pmatrix}$ in the second case.
 \end{enumerate}
 Then we have the list of noncommutative Hom-associative algebras. The commutative case can be established in the same way. In this case the Hom-Ass identity is reduced to
 $$(e_1e_1)f(e_2)-(e_1e_2)f(e_1)=0, \ (e_1e_2)f(e_2)-(e_2e_2)f(e_1)=0.$$
 Then $f$ is in the kernel of the linear system whose matrix is
 $$
 HA_A=
 \begin{pmatrix}
    -\alpha_2\alpha_1-\beta_2\alpha_2  &   -\alpha_2^2 -\beta_2\alpha_4 & \alpha_1^2 + \beta_1\alpha_2 & \alpha_1\alpha_2+ \beta_1\alpha_4\\
    -\alpha_2\beta_1-\beta_2^2  & \alpha_2\beta_2-\beta_2\beta_4 & \alpha_1\beta_1+\beta_1\beta_2 & \alpha_1\beta_2+\beta_1\beta_4 \\
    -\alpha_4\alpha_1-\beta_4\alpha_2 & -\alpha_4\alpha_2-\beta_4\alpha_4 & \alpha_2\alpha_1+\beta_2\alpha _2& \alpha_2^2+\beta_2\alpha_4 \\
    -\alpha_4\beta_1-\beta_4\beta_2 & -\alpha_4\beta_2-\beta_4^2 & \alpha_2\beta_1+\beta_2^2& \alpha_2\beta_2+\beta_2\beta_4  
\end{pmatrix}
$$
Then $A$ is a Hom-associative algebra if and only if $H(A)=\det (HA_A)=0$. We deduce that the set of $2$-dimensional commutative Hom-associative algebra can be provided with an algebraic hypersurface embedded in the affine variety $\K^6$. From Theorem 6, when $\K$ is algebraically closed, we obtain:
\begin{enumerate}
  \item $H(A^6)=\alpha_2\beta_2(1-\alpha_2-\beta_2-3\alpha_2\beta_2+2\alpha_2^2\beta_2+\alpha_2^3\beta_2+2\alpha_2\beta_2^2+2\alpha_2^2\beta^2+\alpha_2\beta_2^3)$. It is equal to $0$ for $\alpha_2=0$ or $\beta_2=0$ or  $\alpha_2=1-\beta_2$ or $\alpha_2=\ds\frac{-3\beta_2-\beta_2^2-(1+\beta_2)\sqrt{\beta_2}\sqrt{4+\beta_2}}{2\beta_2}$ or  $\alpha_2=\ds\frac{-3\beta_2-\beta_2^2+(1+\beta_2)\sqrt{\beta_2}\sqrt{4+\beta_2}}{2\beta_2}.$
  \item $H(A^7)=-\ds\frac{1}{4}$ and $A^7$ is not a Hom-associative algebra.
  \item $H(A^8)=-\ds\frac{9}{64}$ and $A^8$ is not a Hom-associative algebra.
  \item $H(A^i)=0$  for $i=9,10,11,12,13,14,15,16,17$ and $A^9,A^{10},A^{11},A^{12},A^{13},A^{14},A^{15},$ $A^{16},A^{17}$ are  a Hom-associative algebras.
\end{enumerate}


\def\cprime{$'$}

\end{document}